\documentclass[11pt]{article}
\usepackage{amsmath,amsthm,amssymb,xypic}

\theoremstyle{plain}  

\newtheorem{thm}{Theorem}[section]

\newtheorem{lem}[thm]{Lemma}

\theoremstyle{definition}

\newtheorem{prop-defn}[thm]{Proposition--Definition}
\newtheorem{defn}[thm]{Definition}
\newtheorem{notn}[thm]{Notation}

\theoremstyle{remark}

\newtheorem{rem}[thm]{Remark}
\newtheorem{ex}[thm]{Example}

\DeclareMathOperator{\Cl}{Cl}

\DeclareMathOperator{\Spec}{Spec}
\DeclareMathOperator{\Proj}{Proj}

\DeclareMathOperator{\Bl}{Bl}

\DeclareMathOperator{\Supp}{Supp}
\DeclareMathOperator{\Exc}{Ex}

\DeclareMathOperator{\nd}{\not \: \mid}

\newcommand{\QED}{\ifhmode\unskip\nobreak\fi\quad {\rm Q.E.D.}} 

\newcommand{\bA}{\mathbb A}
\newcommand{\bC}{\mathbb C}

\newcommand{\bN}{\mathbb N}
\newcommand{\bP}{\mathbb P}
\newcommand{\bQ}{\mathbb Q}

\newcommand{\bZ}{\mathbb Z}

\newcommand{\cX}{\mathcal X}
\newcommand{\cY}{\mathcal Y}

\newcommand{\cO}{\mathcal O}

\newcommand{\map}{\rightarrow}
\newcommand{\da}{\downarrow}
\newcommand{\inj}{\hookrightarrow}

\title{Semistable divisorial contractions}
\author{Paul Hacking}
\date{June 19th, 2003}
\begin{document}
\maketitle
\pagestyle{plain}

\begin{abstract}
The semistable minimal model program is a special case of the minimal model program concerning 
3-folds fibred over a curve and birational morphisms preserving this structure.
We classify semistable divisorial contractions 
which contract the exceptional divisor to a normal point of a fibre.
Our results can be applied to describe compact moduli spaces of surfaces.

\end{abstract}

\section{Introduction}

The minimal model program is a generalisation of the  classical theory of minimal models of surfaces.
Given a 3-fold $X$ with mild singularities, a minimal model $Y$ of $X$ is a 3-fold birational to $X$
such that either \mbox{$K_Y \cdot C \ge 0$} for every curve $C \subset Y$, or there exists a fibration $Y \map S$ 
such that \mbox{$K_Y \cdot C < 0$} for every curve $C$ contained in a fibre.
The minimal model program constructs a minimal model of $X$ inductively:
if $X$ is not minimal then there is an elementary birational map $X \dashrightarrow X'$ and we
replace $X$ with $X'$. After a finite number of such steps we obtain a minimal model.
The elementary birational maps are of two types: divisorial contractions and flips. 
A divisorial contraction is a birational morphism $\phi \colon X \map Y$ which contracts an
irreducible divisor $E \subset X$ to a point or curve on $Y$. A flip is a birational map 
$\phi \colon X \dashrightarrow Y$ of the form $\phi = g^{-1} \circ f$, where
$f \colon X \map Z$ and $g \colon Y \map Z$ are birational morphisms which contract bunches of
curves to points of $Z$.

Suppose given a 3-fold $X$ fibred over a curve $T$ such that each fibre is reduced
and the total space and the fibres have appropriately mild singularities 
(including, e.g., the case that $X$ is smooth and the fibres are simple normal crossing divisors).
We can run a relative minimal model program for the family $X/T$, contracting only curves which lie in the
fibres, so that each elementary birational map is defined over $T$. This process is called the semistable
minimal model program (\cite{KM}, Chapter~7); we say a contraction is semistable if it occurs in this context.
In this paper we classify semistable divisorial contractions which contract the exceptional divisor to 
a normal point of a fibre.

The semistable minimal model program is an important tool in the construction and explicit description of 
compact moduli spaces of surfaces. For, given a family $\cX^{\times}/T^{\times}$ of surfaces of general type
over a punctured curve $T=T^{\times} \backslash \{ 0 \}$, we can construct a canonical completion to a family $\cX/T$
using the semistable minimal model program \cite{KSB}.
If $M$ is a moduli space  of surfaces of general type, there is a compactification $\bar{M}$ of $M$
with boundary points corresponding to surfaces obtained as the limit of a family $\cX^{\times}/T^{\times}$ 
of surfaces from $M$ by the above process.
More generally, one can compactify moduli spaces of surface-divisor 
pairs $(X,D)$ such that $K_X+D$ is ample \cite{Al}. 
For example, the surface may be a K3 or abelian surface.
In this case we construct limits using the log minimal model program, where the r\^{o}le of 
$K_X$ is played by $K_X+D$.
In fact, I have already applied the results of this paper to give an explicit description of the degenerate 
surfaces which occur at the boundary of a compactification of a moduli space $M$ of pairs (\cite{H1},\cite{H2}).
Here $M$ is the moduli space  of pairs consisting of the plane together with a smooth curve 
of fixed degree $d \ge 4$; in other words $M$ is the moduli space of smooth plane curves of degree $d$.
Similiarly, one can apply our methods to explain and extend the results of Hassett on stable reduction
of plane curve singularities \cite{Has}.

Since this article first appeared as a preprint \cite{H3}, the classification of all divisorial contractions
with centre a point has been completed (\cite{Ha1}, \cite{Ha2}, \cite{K1}, \cite{K2}, \cite{K3}, \cite{K4}).
However, the classification of the semistable contractions is much simpler --- 
there are several exceptional cases which do not occur in the semistable context, 
and the remaining cases are organised into easily described families; moreover, the proof of the classification 
is concise. 
Hence, in view of the applications to moduli problems, it seems worthwhile to present our result separately.

I would like to thank Alessio Corti for his helpful comments on a preliminary version of this paper.

\section{The Classification}

\begin{notn}
We refer to \cite{KM} for Mori theory background, including the definitions of the various classes of 
singularities we shall consider. In this paper log terminal means purely log terminal. 
We write $0 \in T$ for the analytic germ of a smooth curve.
We use script letters to denote flat families of surfaces over $T$
and regular letters for the special fibre, e.g.,
\begin{eqnarray*}
\begin{array}{ccc}
X & \subset & \cX \\
\da &       & \da \\
0   & \in   & T
\end{array}
\end{eqnarray*}
\end{notn}

\begin{defn} 
We say $\cX/T$ is \emph{semistable terminal} if $(\cX,X)$ is
divisorial log terminal and $\cX$ is terminal.
\end{defn}

\begin{defn}
A \emph{semistable divisorial contraction} is a birational morphism 
$\pi \colon \cY \map \cX/T$ where
\begin{enumerate}
\item $\cY/T$ and $\cX/T$ are semistable terminal.
\item The divisor $-K_{\cY}$ is $\pi$-ample.
\item The relative Picard number $\rho(\cY/\cX)=1$.
\item The exceptional locus of $\pi$ is a divisor.
\end{enumerate}
\end{defn}

We recall the classification of semistable terminal singularities (cf. \cite{KSB}). It can also be easily derived using 
the method of toric blowups explained in \cite{YPG}.

\begin{thm} \label{thm-sings}
Let $\cX/T$ be semistable terminal. Then, locally analytically  at a point
$P \in X \subset \cX$,
the pair $X \subset \cX$ is isomorphic to one of the following:
\begin{enumerate}
\item $(xyz=0) \subset \bC^3$.
\item $(t=0) \subset ((xy+tg(z^n,t)=0) \subset \frac{1}{n}(1,-1,a,0))$, 
where $(a,n)=1$. 
\item $(t=0) \subset ((xy-z^{kn}+tg(z^n,t)=0) \subset \frac{1}{n}(1,-1,a,0))$,
where $(a,n)=1$.
\item $(t=0) \subset ((f(x,y,z)+tg(x,y,z,t)=0) \subset \bC^4)$, 
where $((f=0) \subset \bC^3)$ is a Du Val singularity of type $D$ or $E$.
\end{enumerate}
Conversely, any isolated singularity of this form is semistable terminal.
\end{thm}

Note that the special fibre $X$ is normal only in cases (3) and (4). 
In case (3) the singularity $P \in X$ is a cyclic quotient
singularity:
$$X =((xy-z^{kn}=0) \subset \frac{1}{n}(1,-1,a)) \cong \bA^2_{u,v}/\frac{1}{kn^2}(1,kna-1)$$
where $x=u^{kn}$, $y=v^{kn}$ and $z=uv$. If $r=1$ then $X$ is a Du Val singularity of type $A$.
In case (4) the singularity $P \in X$ is a Du Val singularity of type $D$ or $E$.

We can now state our result:
\begin{thm} \label{mainthm}
Let $E \subset \cY \stackrel{\pi}{\map} P \in \cX/T$ be a semistable divisorial contraction such that the special fibre of $\cX/T$ 
is normal.
Then there is a local analytic isomorphism 
$$P \in \cX \cong (f(x,y,z)+tg(x,y,z,t)=0) \subset \frac{1}{n}(1,-1,a,0)$$
such that $\pi$ is given by the weighted blowup of $x,y,z,t$ with weights $w=(w_0,1)$,
where $w(tg) \ge w(f)$. Here $t$ corresponds to a local parameter on the curve $T$, the function $f+tg$ is $\mu_n$-invariant 
and $(a,n)=1$.
We have the following possibilities for $n$, $f$ and $w_0$:
\begin{enumerate}

\item $n \in \bN$ arbitrary:
\begin{enumerate}
\item[$(T)$] $f(x,y,z)=xy-z^{kn}$ and $w_0$ is a primitive vector in the lattice 
$\bZ^3 + \bZ \frac{1}{n}(1,-1,a)$ such that $f$ is homogeneous with respect to the weights $w_0$
of $x,y,z$.
\end{enumerate}
\item $n=1$:
\begin{enumerate}
\item[$(D_m)$] $f(x,y,z)=x^2 + y^2z + z^{m-1}$ and $w_0=(m-1,m-2,2)$, some $m \ge 4$.
\item[$(E_6)$] $f(x,y,z)=x^2 + y^3 + z^4$ and $w_0=(6,4,3)$.
\item[$(E_7)$] $f(x,y,z)=x^2 + y^3 + yz^3$ and $w_0=(9,6,4)$.
\item[$(E_8)$] $f(x,y,z)=x^2 + y^3 + z^5$ and $w_0=(15,10,6)$.
\end{enumerate}
\end{enumerate}
Conversely, let
\mbox{$P \in \cX \cong (f+tg=0) \subset \frac{1}{n}(1,-1,a,0)$} and let $\pi\colon \cY \map \cX$ be the 
birational map induced
by the weighted blowup of $x,y,z,t$ with weights $w=(w_0,1)$, where $n$, $f$ and $w_0$ are as above and $g$ is chosen so that 
$\cX$ has an isolated singularity and $w(tg) \ge w(f)$. Then the families $\cY/T$ and $\cX/T$ are 
semistable terminal, the divisor $-K_{\cY}$ is $\pi$-ample and the exceptional locus of $\pi$ is an irreducible divisor.
Thus $\pi$ is a semistable divisorial contraction if \mbox{$\rho(\cY/\cX)=1$} (which is automatic if 
$\cX$ is $\bQ$-factorial at $P$).  
\end{thm}

\section{Proof of the Classification}

\begin{prop-defn}
Let $\pi \colon \cY \map \cX/T$ be a semistable divisorial contraction with 
exceptional divisor $E$.
Assume that $E$ is contracted to a point $P \in \cX$ and that $X$ is normal.
Write $Y=Y_1+E$, where $Y_1$ is the strict transform of $X$ and $F=E |_{Y_1}$.
Then the map $p \colon Y_1 \map X$ satisfies the following:
\begin{enumerate}
\item The surface $X$ is log terminal and the pair $(Y_1,F)$ is log terminal. 
\item The divisor $-(K_{Y_1}+F)$ is $p$-ample.
\item The relative Picard number $\rho(Y_1/X)$ equals $1$.
\item The exceptional locus of $p$ is $F$, a divisor.
\end{enumerate}
We say a map $p \colon Y_1 \map X$ satisfying the conditions above is a 
\emph{log divisorial contraction}.
\end{prop-defn}
\begin{proof}
The pair $(\cX,X)$ is dlt and  $X$ is normal by assumption, 
hence $X$ is log terminal by adjunction.
We have $(K_{\cY}+Y)|_{Y_1}=K_{Y_1}+F$ by adjunction.
The pair $(\cY,Y)$ is dlt, so  $F$ is smooth and irreducible
and $(Y_1,F)$ is log terminal.
The divisor $-K_{\cY}$ is $\pi$-ample thus the restriction
$-(K_{Y_1}+F)$ is $p$-ample. 
\end{proof}

\begin{thm} \label{log_div_contr}
Let $P \in X$ be a log terminal surface singularity which admits a 
$\bQ$-Gorenstein smoothing (equivalently, $X$ occurs as the special fibre of a
semistable terminal family $\cX/T$). Then the log divisorial contractions $E \subset Y_1 \map P \in X$
are precisely the following:
\begin{enumerate}
\item There is an isomorphism $P \in X \cong \bA^2_{u,v} / \frac{1}{kn^2}(1,kna-1)$, where $(a,n)=1$, 
such that
$Y_1 \map X$ is given by the weighted blowup of $u,v$ with weights $\alpha$. 
Here $\alpha \in N = \bZ^2+\bZ\frac{1}{kn^2}(1,kna-1)$ is a primitive vector in the 
lattice $N$.

Equivalently, there is an isomorphism $P \in X \cong (xy-z^{kn}=0) \subset \frac{1}{n}(1,-1,a)$
such that $Y_1 \map X$ is given by the weighted blowup of $x,y,z$ with weights $w_0$.
Here $w_0$ is a primitive vector in the lattice 
\mbox{$\bZ^3 + \bZ \frac{1}{n}(1,-1,a)$} such that $f$ is homogeneous with respect to these weights.

\item The singularity $P \in X$ is a Du Val singularity of type $D$ or $E$.
Then $P \in X \cong ((f(x,y,z)=0) \subset \bC^3)$ and 
$Y_1 \map X$ is given by the weighted blowup of $x,y,z$ with weights $w_0$, where
\begin{enumerate}
\item[$(D_m)$] $f(x,y,z)=x^2 + y^2z + z^{m-1}$ and $w_0=(m-1,m-2,2)$, some $m \ge 4$.
\item[$(E_6)$] $f(x,y,z)=x^2 + y^3 + z^4$ and $w_0=(6,4,3)$.
\item[$(E_7)$] $f(x,y,z)=x^2 + y^3 + yz^3$ and $w_0=(9,6,4)$.
\item[$(E_8)$] $f(x,y,z)=x^2 + y^3 + z^5$ and $w_0=(15,10,6)$.
\end{enumerate}
\end{enumerate}
\end{thm}

\begin{proof}
The pair $(Y_1,F)$ has log terminal singularities, thus at $F$ the singularities of $(Y_1,F)$ 
are of the form $(\frac{1}{r}(1,a),(x=0))$ (\cite{KM}, p. 119, Theorem~4.15).
Let $\tilde{X} \map X$ and $\tilde{Y_1} \map Y_1$ be the minimal resolutions of $X$ and $Y_1$, then there is a diagram
\begin{eqnarray*}
\begin{array}{ccc}
\tilde{Y_1} & \map & \tilde{X}  \\
\da &       & \da \\
Y_1   & \map   & X .
\end{array}
\end{eqnarray*}
The map $q \colon \tilde{Y_1} \map \tilde{X}$ is a composition of contractions of $(-1)$-curves.
The exceptional locus of $\tilde{Y_1} \map X$ consists of the strict transform $F'$ of $F$ 
together with strings of curves meeting $F'$
obtained by resolution of the \mbox{$(\frac{1}{r}(1,a), (x=0))$}
singularities of $(Y_1,F)$. Thus if $q$ is not an isomorphism, 
then $F'$ must be a $(-1)$-curve, and $q$ factors through 
the blow down \mbox{$\sigma \colon \tilde{Y_1} \map Y_1'$} of $F'$.

The exceptional locus of $\tilde{X} \map X$ is either a string of rational curves, or one of the $D$ or $E$  configurations,
in particular it has normal crossing singularities. This property is preserved by blowups, thus if $F'$ is contracted by $q$
then the exceptional locus of $Y_1' \map X$ has normal crossing singularities --- equivalently, the exceptional locus
of $\tilde{Y_1} \map X$ is a string of curves. In this case, the singularity $P \in X$ is a cyclic quotient singularity, and
is 
of the form $\bA^2_{u,v} / \frac{1}{kn^2}(1,kna-1)$ by the classification of Theorem~\ref{thm-sings}. 
Moreover, for an appropriate choice of the coordinates $u,v$, the morphism
$\tilde{Y_1} \map X$ is toric. Then in particular $Y_1 \map X$ is toric and thus is given by a weighted blowup of $u,v$ with 
weights given by some primitive vector in the lattice $\bZ^2+\bZ\frac{1}{kn^2}(1,kna-1)$. 

On the other hand, if the exceptional locus of $\tilde{Y_1} \map X$ is not a string of curves, then 
$q \colon \tilde{Y_1} \map \tilde{X}$ is an isomorphism and the curve $F'$ is the `fork' curve of the $D$ or $E$ 
configuration (i.e., the exceptional curve which meets $3$ others). 
This determines $Y_1 \map X$ uniquely, and we observe that these contractions can be described explicitly as above.
\end{proof}
  
Our main theorem~\ref{mainthm} follows immediately from Theorem~\ref{thm-inversion} below.

\begin{thm} \label{thm-inversion}
Let $E \subset \cY \stackrel{\pi}{\map} P \in \cX/T$ be a semistable divisorial contraction with $X$ normal.
Let $F \subset Y_1 \stackrel{p}{\map} P \in X$ be the induced log divisorial contraction.
There is an isomorphism $$P \in X \cong (f(x,y,z)=0) \subset \frac{1}{n}(1,-1,a)$$  
such that 
$p$ is given by the  weighted blowup of $x,y,z$ with weights $w_0$
as in Theorem~\ref{log_div_contr} --- fix one such identification.
Then there is a compatible isomorphism $$P \in \cX/T \cong ((f(x,y,z)+tg(x,y,z,t)=0) \subset \frac{1}{n}(1,-1,a,0))/\bC^1_t$$
such that $\pi$ is given by the weighted blowup of $x,y,z,t$ with weights $w=(w_0,1)$,
where $w(tg) \ge w(f)$.

Conversely, let $\pi \colon \cY \map \cX$ be a birational morphism constructed in this fashion, 
then the family $\cY/T$ is 
semistable terminal, the divisor $-K_{\cY}$ is $\pi$-ample and the exceptional locus of $\pi$ is an irreducible
 divisor. 
More precisely, we assume that
\mbox{$P \in \cX \cong (f+tg=0) \subset \frac{1}{n}(1,-1,a,0)$} and $\pi$ is given by the weighted blowup of 
$x,y,z,t$ with weights
$w=(w_0,1)$, where \mbox{$(f(x,y,z)=0)\subset\frac{1}{n}(1,-1,a)$} and $w_0$ are as in
Theorem~\ref{log_div_contr}, and $g$ is chosen so that $\cX$ has an isolated singularity and
$w(tg) \ge w(f)$.
\end{thm}

\begin{proof}
A divisorial contraction $\pi \colon E \subset Y \map P \in X$ determines a discrete valuation
$\nu \colon k(X)^{\times} \map \bZ$, where $\nu(f)$ equals the order of vanishing of $\pi^{\star}f$ along $E$.
Moreover, we can reconstruct $\pi$ from $\nu$. For, writing
$$m^{(n)}_{X,P}=\{f \in \cO_{X,P} \mbox{ $|$ } \nu(f) \ge n \},$$
we have $m^{(n)}_{X,P}= \pi_{\star}\cO_Y(-nE)$ and 
$$Y=\underline{\Proj}_X(\oplus_{n \ge 0} \, \pi_{\star}\cO_Y(-nE))$$
since $-E$ is $\pi$-ample. Hence, in order to identify two divisorial contractions it is enough to identify the corresponding
valuations.

Consider the divisorial contraction $p \colon Y_1 \map X$; let $\nu_0 \colon k(X)^{\times} \map \bZ$ be the corresponding 
valuation. Write $w_0=\frac{1}{d}(a_1,a_2,a_3)$, where $a_1$, $a_2$, and $a_3$ are coprime.
We have an inclusion 
$$X=(f(x,y,z)=0) \subset \bA=\frac{1}{n}(1,-1,a)$$
and a weighted blowup $s \colon \Bl_{\frac{1}{d}(a_1,a_2,a_3)}\bA \map \bA$ inducing the contraction $p$.
Write $w_0 \colon k(\bA)^{\times} \map \bZ$ for the valuation corresponding to the weighted blowup;
for $h = \sum a_{ijk}x^iy^jz^k  \in \cO_{\bA,0}$ we have 
$$w_0(h)= \min \{ \frac{1}{d}(a_1 i+a_2 j +a_3 k) \mbox{ $|$ }a_{ijk} \neq 0 \}.$$ 
Then, for $h \in \cO_{X,P}$, 
$$\nu_0(h)= \max\{w_0(\tilde{h}) \mbox{ $|$ } \tilde{h} \in \cO_{\bA,0} \mbox{ a lift of $h$}\}.$$ 
More precisely, given a lift $\tilde{h}$ of $h$, write 
$$\tilde{h}=\tilde{h}_{\alpha}+\tilde{h}_{\alpha+1}+\cdots,\,\tilde{h}_{\alpha} \neq 0,$$
for the decomposition of $\tilde{h}$ into graded pieces with respect to the weighting.
Then $\nu_0(h)=w_0(\tilde{h})$ iff $f \nd \tilde{h}_{\alpha}$.
For, the exceptional divisor $G$ of $s$ is a weighted projective space $\bP(a_1,a_2,a_3)$,
and the exceptional divisor $F$ of $p$ is given by 
$$F=(f(X,Y,Z)=0) \subset  \bP(a_1,a_2,a_3)$$
(note that $f$ is homogeneous with respect to the weighting).
Now, $s^{\star}\tilde{h}=u^{\alpha}\tilde{h}'$, where $u$ is a local parameter at $G$, and
$(\tilde{h}'=0)$ is the strict transform of $(\tilde{h}=0) \subset \bA$.
Thus $p^{\star}h=s^{\star}\tilde{h}|_{Y_1}=v^{\alpha} \cdot \tilde{h}'|_{Y_1}$,
where $v$ is a local parameter at $F$. So $\nu_0(h) \ge \alpha=w_0(\tilde{h})$, with equality
iff $\tilde{h}'|_{Y_1}$ does not vanish on $F$. We have
$$(\tilde{h}'=0)|_G=(\tilde{h}_{\alpha}(X,Y,Z)=0) \subset \bP(a_1,a_2,a_3),$$
so this last condition is equivalent to $f \nd \tilde{h}_{\alpha}$.

Our aim is to deduce a description of the contraction $\pi \colon \cY \map \cX$ similiar to that of
$p \colon Y_1 \map X$ above. Let $\nu$ be the valuation defined by $\pi$. We first show that, in the case
that the index $n$ of $P \in \cX$ equals $1$, we may lift $x,y,z \in \cO_{X,P}$ to 
$\tilde{x},\tilde{y},\tilde{z} \in \cO_{\cX,P}$ such that $\nu(\tilde{x})=\nu_0(x)$ etc.
Certainly, for $h \in \cO_{X,P}$ and any lift $\tilde{h} \in \cO_{\cX,P}$, we have $\nu(\tilde{h}) \le \nu_0(h)$;
we show that for an appropriate choice of $\tilde{h}$ we have equality. Equivalently, the map 
$$\pi_{\star}\cO_{\cY}(-iE) \map p_{\star}\cO_{Y_1}(-iE)$$
is surjective for each $i \ge 1$. Applying $\pi_{\star}$ to the exact sequence
$$0 \map \cO_{\cY}(-iE-Y_1) \map \cO_{\cY}(-iE) \map \cO_{Y_1}(-iE) \map 0$$
we obtain the long exact sequence
$$\cdots \map \pi_{\star}\cO_{\cY}(-iE) \map p_{\star}\cO_{Y_1}(-iE) \map R^1 \pi_{\star}\cO_{\cY}(-iE-Y_1) \map \cdots.$$
We have $Y_1+E = Y \sim 0$ and $K_{\cY} \sim \pi^{\star}K_{\cX}+aE \sim aE$, so $-iE-Y_1 \sim K_{\cY}-(i-1+a)E$. 
Thus $R^1\pi_{\star}\cO_{\cY}(-iE-Y_1)=R^1\pi_{\star}\cO_{\cY}(K_{\cY}-(i-1+a)E)=0$ by Kodaira vanishing since $-E$ is
$\pi$-ample; the required surjectivity result follows.
In the case that the index is greater than $1$, by Lemma~\ref{covering_lem} there is a diagram
\begin{eqnarray*}
\begin{array}{ccc}
\tilde{E} \subset \tilde{\cY} & \stackrel{\tilde{\pi}}{\map} & \tilde{P} \in \tilde{\cX}\\
\da                             &                            & \da\\
E \subset \cY &  \stackrel{\pi}{\map}                        & P \in \cX\\
\end{array}
\end{eqnarray*} 
where $\tilde{\cX} \map \cX$ is the index one cover. Then a similiar calculation shows that we may lift
$x,y,z \in \cO_{\tilde{X},P}$ to $\tilde{x},\tilde{y},\tilde{z} \in \cO_{\tilde{\cX},P}$ of the same weights
as above.

For simplicity of notation, now write $x,y,z$ for the lifts $\tilde{x},\tilde{y},\tilde{z} \in \cO_{\tilde{\cX},P}$ of 
$x,y,z \in \cO_{\tilde{X},P}$ constructed above. Let $t$ denote a local parameter at $0 \in T$. Then $x,y,z,t$ define an 
embedding $\cX \inj \bA \times T$ extending $X \inj \bA$; write
$$\cX=(f(x,y,z)+tg(x,y,z,t)=0) \subset \bA \times T = \frac{1}{n}(1,-1,a,0).$$
Define a weighted blowup 
$$\sigma \colon \Bl_{\frac{1}{d}(a_1,a_2,a_3,d)}(\bA \times T) \map \bA \times T,$$
and let $w$ denote the corresponding valuation. Our aim is to show that $\pi$ is induced by $\sigma$.
Note immediately that, for $h \in \cO_{X,P}$, and any lift $\tilde{h} \in \cO_{\bA \times T,0}$, we have
$\nu(h) \ge w(\tilde{h})$.
For, the valuations $\nu$ and $w$ agree on $x,y$ and $z$ by construction, and $\nu(t)=1=w(t)$ by the semistability assumption.
So, writing $\tilde{h}=\sum a_{ijkl}x^iy^jz^kt^l$, we have
$$\nu(h) \ge \min \{ \nu(x^iy^jz^kt^l) \mbox{ $|$ } a_{ijkl} \neq 0 \} = \min \{ w(x^iy^jz^kt^l) \mbox{ $|$ } a_{ijkl} \neq 0 \} = w(\tilde{h})$$
as claimed. We also need the following preliminary result: $w(tg) \ge w(f)$.
To prove this, write 
$$g=g_{\alpha}+g_{\alpha+1}+ \cdots,\, g_{\alpha} \neq 0$$ 
for the decomposition of $g$ into graded pieces with respect to the weighting. Let $g_{\alpha}=t^{\beta}k$, where
$k_0=k|_{t=0} \neq 0$.  
We may assume that $f$ does not divide $k_0$, for otherwise $w(tg)>w(f)$.
In this case, we have $\nu_0(k_0)=w_0(k_0)$ as proved above.
Now $w(k) \le \nu(k) \le \nu_0(k_0) = w_0(k_0) =w(k)$, thus $w(k)=\nu(k)$ and so $w(g)=\nu(g)$.
Finally $w(f) \le \nu(f)=\nu(tg)=w(tg)$ as required.

We now complete the proof that $\pi$ is induced by $\sigma$. The contraction $\pi' \colon E' \subset \cY' \map P \in \cX$
induced by $\sigma$ has valuation $\nu'$ given by
$$\nu'(h)= \max\{w(\tilde{h}) \mbox{ $|$ } \tilde{h} \in \cO_{\bA \times T,0} \mbox{ a lift of $h$}\}$$
(cf. our earlier treatment of the contraction $p$ induced by $s$).
We show that $\nu=\nu'$ and thus $\pi=\pi'$ as required.
We know that $\nu(h) \ge w(\tilde{h})$ for $h \in \cO_{X,P}$ and $\tilde{h} \in \cO_{\bA\times T, 0}$ any lift of $h$,
it remains to show that we have equality for some $\tilde{h}$.
Given $h \in \cO_{\cX,P}$, pick some lifting $\tilde{h}$. Write
$$\tilde{h} =\tilde{h}_{\alpha}+\tilde{h}_{\alpha+1}+\cdots$$
for the decomposition into graded pieces with respect to the weighting and 
$\tilde{h}_{\alpha}=t^{\beta}k$, where $k_0=k|_{t=0} \neq 0$.
If $f$ divides $k_0$, say $k_0=q \cdot f$, replace $\tilde{h}$ by $\tilde{h}'=\tilde{h}-q(f+tg)$.
Defining $\alpha'$ and $\beta'$ as above, we see that either $\alpha' > \alpha$
or $\alpha'=\alpha$ and $\beta' > \beta$, using $w(tg) \ge w(f)$. 
It follows that this process can repeat only finitely many times, so we may assume
$f \nd k_0$ and hence $\nu_0(k_0)=w_0(k_0)$. 
Then $w(k)\le \nu(k) \le \nu_0(k_0)=w_0(k_0)=w(k)$, so $w(k)=\nu(k)$ and $w(\tilde{h})=\nu(h)$ as required.

Conversely, we show that all contractions $\pi \colon \cY \map \cX$ constructed in this way have all the 
properties of a semistable divisorial contraction, except that the relative Picard number is not necessarily equal to $1$.
We first prove that $(\cY,Y)$ is dlt. We know that the pair $(Y_1,F)$ is log terminal and
$$E=(f(X,Y,Z)+Tg_{\lambda-1}(X,Y,Z,T)=0) \subset \bP(a_1,a_2,a_3,d),$$
where $\lambda=w(f)$ and $g_{\lambda-1}$ is the graded piece of $g$ with weight $\lambda-1$ (possibly zero).
The curve $F=(T=0)\subset E$ is smooth and the only singularities of $E$ at $F$ are cyclic quotient singularities 
induced by the singularities of the ambient weighted projective space. Hence the pair $(E,F)$ is log terminal near $F$.
But $-(K_E+F)$ is ample, so $(E,F)$ is log terminal everywhere by Shokurov's connectedness result (see Lemma~\ref{Shokurov}).
Hence the degenerate surface $Y=Y_1+E$ is semi log terminal \cite{KSB}. 
A standard inversion of adjunction argument shows that $(\cY,Y)$ is dlt.
The divisor $E$ is clearly irreducible. Finally, $-K_{\cY}$ is $\pi$-ample. For $K_{\cY} \equiv aE$ 
and $a>0$ since $\cX$ is terminal. 
\end{proof}

\begin{rem}
Note that a morphism $\pi \colon \cY \map \cX$ constructed as in the Theorem is a semistable divisorial contraction
precisely when $\rho(\cY/\cX)=1$. This is not always the case, e.g., if $\cX=(xy+z^2+t^2=0)\subset \bC^4$
and $\pi \colon \cY \map \cX$ is the blowup of $0 \in \cX$ then $\rho(\cY/\cX)=2$, and $\pi$ can be factored
into a divisorial contraction followed by a flopping contraction. However, if $\cX$ is $\bQ$-factorial, then
the exact sequence
$$0 \map \bZ E \map \Cl(\cY) \map \Cl(\cX) \map 0$$
implies that $\rho(\cY/\cX)=1$. 

Note also that, in order to compute $\rho(\cY/\cX)$, we must work algebraically rather than local analytically at $P \in \cX$.
For example, we can construct a variety $\cX$ which has an ordinary double point $P \in \cX$ but is (algebraically) $\bQ$-factorial. 
Then the blowup $\cY \map \cX$ of $P \in \cX$ has $\rho(\cY/\cX)=1$. 
\end{rem}

\begin{lem}\label{covering_lem}
Let $\pi \colon E \subset \cY \map P \in \cX$ be a semistable divisorial contraction.
Let $p \colon \tilde{\cX} \map \cX$ be the index one cover of $\cX$.
Then there is a diagram
\begin{eqnarray*}
\begin{array}{ccc}
\tilde{E} \subset \tilde{\cY} & \stackrel{\tilde{\pi}}{\map} & \tilde{P} \in \tilde{\cX}\\
q \da                           &                            & p \da\\
E \subset \cY &  \stackrel{\pi}{\map}                        & P \in \cX\\
\end{array}
\end{eqnarray*}
where the map $q \colon \tilde{\cY} \map \cY$ is a cyclic quotient,
and the map $\tilde{\pi} \colon \tilde{\cY} \map \tilde{\cX}$ is a birational morphism 
such that $\tilde{\cY}$ is log terminal, the divisor $-K_{\tilde{\cY}}$ is $\tilde{\pi}$-ample
and the exceptional locus $\Exc(\tilde{\pi})=\tilde{E}$ is an irreducible divisor.
\end{lem}

\begin{proof}
Note first that $\tilde{\cY}$ is the normalisation of $\cY$ in the function field of $\tilde{\cX}$. 
We give an explicit construction of $\tilde{\cY}$ below and verify the desired properties.

Let $a$ be the discrepancy of $E$ and $n$ the index of $P \in \cX$. 
Then the index one cover of $P \in \cX$ is given by
$$\tilde{\cX}=\underline{\Spec}_{\cX}(\oplus_{i=0}^{n-1} \cO_{\cX}(iK_{\cX})),$$
where the multiplication is defined by fixing an isomorphism $\cO_{\cX}(nK_{\cX}) \cong \cO_{\cX}$.
Define $\tilde{\cY}$ as follows:
$$\tilde{\cY}=\underline{\Spec}_{\cY}(\oplus_{i=0}^{n-1} \cO_{\cY}(iK_{\cY}+ \lfloor -iaE \rfloor)),$$
where the multiplication is given by $\cO_{\cY}(nK_{\cY}-naE)=\cO_{\cY}(\pi^{\star} nK_{\cX}) \cong \cO_{\cY}$.
The map $\tilde{\pi} \colon \tilde{\cY} \map \tilde{\cX}$ is given by the natural maps
$$\cO_{\cX}(iK_{\cX}) \map \pi_{\star}\cO_{\cY}(iK_{\cY}+ \lfloor -iaE \rfloor),$$
where we are using $\pi^{\star}K_{\cX} \equiv K_{\cY}-aE$. 

Write $a=a_1/d$ where $a_1$ and $d$ are coprime integers. Then $d$ divides $n$; let $e=n/d$.
The exact sequence of groups
$$0 \map \mu_d \map \mu_n \map \mu_e \map 0$$
induces a factorisation $\tilde{\cY} \map \cY' \map \cY$ where $\cY' \map \cY$ is etale in codimension~$1$: 
we have
$$\cY' =\tilde{\cY}/\mu_d = \underline{\Spec}_{\cY}(\oplus_{j=0}^{e-1} \cO_{\cY}(jdK_{\cY}-ja_1E))$$
and
$$\tilde{\cY}= \underline{\Spec}_{\cY'}(\oplus_{k=0}^{d-1} \cO_{\cY'}(kK_{\cY'}+\lfloor -kaE' \rfloor)).$$
To understand  the cover $\tilde{\cY} \map \cY'$,
choose $b$ such that $-ba_1=1 \mod d$ and, working locally over a point of $E' \subset \cY'$ where 
$\cY'$ is smooth, let $v$ generate $\cO_{\cY'}(bK_{\cY'}+\lfloor -baE' \rfloor)$.
Then $\cO_{\tilde{\cY}}$ is freely generated locally by $1,v,\ldots, v^{d-1}$ over $\cO_{\cY'}$ and $v^d=u \in 
\cO_{\cY'}$ cuts out $E'$. In other words,
$$\tilde{\cY}=\Spec \cO_{\cY'}[s]/(s^d=u)$$
where $E'=(u=0)$. So $\tilde{\cY}$ is normal and $\tilde{\cY} \map \cY'$ is totally ramified over $E'$.

We verify that $\tilde{\cY}$ is log terminal.
The pair $(\cY,Y)$ is dlt by the semistability assumption, in particular $(\cY,E)$ is log terminal.
We have $K_{\tilde{\cY}}=q^{\star}K_{\cY} + (d-1)\tilde{E}$ by Riemann--Hurwitz  and 
$q^{\star}E=d\tilde{E}$, so $K_{\tilde{\cY}}+\tilde{E}=q^{\star}(K_{\cY}+E)$.
Thus $(\cY,E)$  log terminal implies $(\tilde{\cY},\tilde{E})$ log terminal and so $\tilde{\cY}$ is log terminal
as required. The remaining assertions are clear.
\end{proof}

\begin{rem}
The contraction $\tilde{\pi}$ covering $\pi$ is usually not a Mori contraction since $\tilde{\cY}$ has log terminal 
but not terminal singularities. However, in some cases, $\tilde{\pi}$ is a Mori contraction.
For example, let 
$$\cX=(xy+z^2+t^N=0)\subset \frac{1}{2}(1,1,1,0), N \ge 3,$$
a $\mu_2$ quotient of a $cA_1$ point, and let $\pi$ be given by the weighted blowup with weights $\frac{1}{2}(1,5,3,2)$.
Then $\tilde{\cX}=(xy+z^2+t^N=0)\subset\bA^4$
and $\tilde{\pi}$ is the contraction given by the weighted blowup with weights $(1,5,3,2)$. 
A calculation shows that $\tilde{\pi}$ is a Mori contraction precisely when $N=3$. 
Moreover this is the only non-semistable contraction to a $cA_1$ point \cite{K2}.
\end{rem}

\begin{lem} \label{Shokurov}
Let $X$ be a normal proper surface and $B$ an effective $\bQ$-divisor on $X$ such that $-(K_X+B)$ is nef and big.
Then the locus where $(X,B)$ is not klt is connected.
\end{lem}

\begin{proof} 
This is a global version of a connectedness result due to Shokurov and Koll\'{a}r (see e.g. \cite{KM}, p. 173, Theorem~5.48)
--- the proof is essentially unchanged but is included below for completeness.
 
Let $f \colon Y \map X$ be a resolution of $X$ such that the support of the strict transform of $B$ together with 
the exceptional locus of $f$ form a snc divisor on $Y$. Define $D$ via the equation
$$K_Y+D = f^{\star}(K_X+B)$$
and let $D=\sum d_i D_i$ where the $D_i$ are the irreducible components of $\Supp D$.
Decompose $D$ into $A=\sum_{d_i < 1} d_i D_i$ and $F=\sum_{d_i \ge 1} d_i D_i$.
Then the locus where $(X,B)$ is not klt is the image of $F$, so we are required to prove that 
$\Supp F$ is connected.
By definition we have
$$-A-F=-D=K_Y-f^{\star}(K_X+B).$$
Rounding up we obtain
$$\lceil -A \rceil - \lfloor F \rfloor = K_Y - f^{\star}(K_X+B)+ \{A\} + \{F\}.$$  
Since $-f^{\star}(K_X+B)$ is nef and big and $\{A\}+\{F\}$ has snc support and coefficients smaller than $1$,
we deduce that $H^1(\cO_Y(\lceil -A \rceil - \lfloor F \rfloor))=0$ by Kodaira vanishing.
Hence the map
$$H^0(\cO_Y(\lceil -A \rceil)) \map H^0(\cO_{\lfloor F \rfloor}(\lceil -A \rceil))$$
is surjective. The divisor $\lceil -A \rceil$ is effective and $f$-exceptional.
So $$H^0(\cO_Y(\lceil -A \rceil))=H^0(\cO_X) \cong \bC$$ and $H^0(\cO_{F_j}(\lceil -A \rceil)) \neq 0$
for any connected component $F_j$ of $\lfloor F \rfloor$.
Thus $\Supp \lfloor F \rfloor = \Supp F$ is connected by the surjectivity result above, as required.
\end{proof}

\section{Examples}

We give a detailed description of the semistable contractions \mbox{$E \subset \cY \stackrel{\pi}{\map} P \in \cX$} 
where $P \in \cX$ is of the form
$$(xy+z^{kn}+tg(x,y,z,t)=0) \subset \frac{1}{n}(1,-1,a,0).$$
By our theorem we may choose the coordinates $x,y,z,t$ so that the contraction $\pi$ is given by the weighted blowup of
$x,y,z,t$ with weights $w=(w_0,1)$, where $w_0 \in \bZ^3+\bZ\frac{1}{n}(1,-1,a)$ is a primitive vector, 
$f(x,y,z)=xy+z^{kn}$ is homogeneous with respect to $w_0$ and $w(tg) \ge w(f)$.

\begin{ex}
Consider first the case where the index $n$ equals $1$. Let $w_0=(a_1,a_2,a_3)$.
After a change of coordinates we may assume that 
$$tg(x,y,z,t)=b_{k-2}z^{k-2}t^{2a_3}+b_{k-3}z^{k-3}t^{3a_3}+\cdots+b_0t^{ka_3},$$
where $b_i \in k[[t]]$ for each $i$, using $w(tg) \ge w(f)$.
Then the exceptional divisor $E$ of $\pi$ is the surface
$$E=(XY+Z^k+c_{k-2}Z^{k-2}T^{2a_3}+\cdots+c_0T^{ka_3}=0) \subset \bP(a_1,a_2,a_3,1)$$
where $c_i=b_i(0) \in k$, and $X,Y,Z,T$ are the homogeneous coordinates on $\bP(a_1,a_2,a_3,1)$ corresponding to
the coordinates $x,y,z,t$ on $\bA^4$. 

We now describe the singularities of the pair $(\cY,Y)$.
Consider first the affine piece $U=(T\neq 0)$ of $E$:
$$U = (x'y'+{z'}^k+c_{k-2}{z'}^{k-2}+ \cdots + c_0 =0) \subset \bA^3.$$
The only singularities of $U$ are $A_{l-1}$ singularities at the points 
where $x'=y'=0$ and $h(z')={z'}^k+c_{k-2}{z'}^{k-2} + \cdots + c_0$ has a multiple root
of multiplicity $l$. These singularities are sections of $cA_{l-1}$ singularities on the total space $\cY$.
The remaining singularities of $Y \subset \cY$ are of the form $(xy=0) \subset \frac{1}{a_1}(1,-1,a_3)$
and $(xy=0) \subset \frac{1}{a_2}(1,-1,a_3)$ at the points $(1:0:0:0)$ and $(0:1:0:0)$ of $E$ respectively.
\end{ex}

\begin{ex}
Let the index $n \in \bN$ be arbitrary. Let $w_0=\frac{1}{d}(a_1,a_2,a_3)$, where $a_1$, $a_2$, and $a_3$ are coprime,
and write $n=de$. We may assume that
$$tg(x,y,z,t)=b_{k-2}z^{(k-2)n}t^{2ea_3}+b_{k-3}z^{(k-3)n}t^{3ea_3}+\cdots+b_0t^{kea_3},$$
where $b_i \in k[[t]]$ for each $i$, and we have used the fact that $g$ is $\mu_n$ invariant.
The exceptional divisor $E$ of $\pi$ is the surface
$$E=(XY+Z^{kn}+c_{k-2}Z^{(k-2)n}T^{2ea_3}+\cdots+c_0T^{kea_3}=0) \subset \bP(a_1,a_2,a_3,d)$$
where $c_i=b_i(0) \in k$. 
Consider the affine piece $U=(T\neq 0)$ of $E$:
$$U = (x'y'+{z'}^{kn}+c_{k-2}{z'}^{(k-2)n}+ \cdots + c_0 =0) \subset \frac{1}{d}(a_1,a_2,a_3).$$
The singularities of $U$ away from $(0,0,0)$ can be described as above. At $(0,0,0)$ we have a singularity
of the form
$$(xy+z^{ln}=0) \subset \frac{1}{d}(a_1,a_2,a_3) \cong \frac{1}{d}(1,-1,b)$$
where $l$ is the least $i$ such that $b_i \neq 0$ and $b \equiv a_1^{-1}a_3 \mod d$.
This is a cyclic quotient singularity of the form $\frac{1}{k'{n'}^2}(1,k'n'a'-1)$, where $k'=le$, $n'=d$  and $a'=b$.
The corresponding singularity of the total space $\cY$ is a deformation of $0 \in U$ of the form 
$$(xy+z^{ln}+tg(z^d,t)=0) \subset \frac{1}{d}(1,-1,b,0).$$
The remaining singularities of $Y \subset \cY$ are of the form $(xy=0) \subset \frac{1}{ea_1}(1,-1,\frac{a_3-aa_1}{d})$
and $(xy=0) \subset \frac{1}{ea_2}(1,-1,\frac{aa_2+a_3}{d})$ at the points $(1:0:0:0)$ and $(0:1:0:0)$ of $E$ respectively.
\end{ex}

\noindent Department of Mathematics,\\
University of Michigan,\\ 
Ann Arbor, MI 48109, USA\\ 
phacking@umich.edu\\


\begin{thebibliography}{99}

\bibitem[1]{Al} V. Alexeev, Moduli spaces $M_{g,n}(W)$ for surfaces, in 
Higher dimensional complex varieties, Trento, 1994, 1--22.
\bibitem[2]{H1} P. Hacking, A compactification of the space of plane curves, Cambridge Univ. PhD thesis, 
2001, math.AG/0104193. 
\bibitem[3]{H2} P. Hacking, Compact moduli of plane curves, submitted. 
\bibitem[4]{H3} P. Hacking, Semistable divisorial contractions, preprint, math.AG/0208049.
\bibitem[5]{Has} B. Hassett, Local stable reduction of plane curve singularities, J. Reine Angew. Math. 520 (2000), 169--194
\bibitem[6]{Ha1} T. Hayakawa, Blowing ups of 3 dimensional terminal singularities, Publ. Res. Inst. Math. Sci. 35 (1999), 515--570.
\bibitem[7]{Ha2} T. Hayakawa, Blowing ups of 3 dimensional terminal singularities II, Publ. Res. Inst. Math. Sci. 36 (2000), 423--456.
\bibitem[8]{K1} M. Kawakita, Divisorial contractions in dimension $3$ which contract divisors to smooth 
points, Invent. Math. 145 (2001), no. 1, 105--119.
\bibitem[9]{K2} M. Kawakita, Divisorial contractions in dimension $3$ which contract divisors to 
compound $A_1$ points, Compositio Math. 133 (2002), 95--116.
\bibitem[10]{K3} M. Kawakita, General elephants of three-fold divisorial contractions, J. Amer. Math. Soc. 16 
(2003), 331--362.
\bibitem[11]{K4} M. Kawakita, Three-fold divisorial contractions to singularities of higher indices, 
preprint, math.AG/0306065. 
\bibitem[12]{KM} J. Koll\'{a}r and S. Mori, Birational geometry of algebraic varieties, C.U.P., 1998.
\bibitem[13]{KSB} J. Koll\'{a}r and N. Shepherd-Barron, Threefolds and deformations of surface singularities, Invent. Math. 91 (1988), no. 2, 299--338.
\bibitem[14]{YPG} Young person's guide to canonical singularities, in Algebraic Geometry, Bowdoin 1985, ed. S. Bloch, Proc. of Symposia in Pure Math. 46, A.M.S. (1987), vol. 1, 345--414.
\end{thebibliography}
\end{document}